\documentclass[12pt,reqno]{amsart}
\swapnumbers

\theoremstyle{plain}
\newtheorem*{theorem}{Theorem}
\newtheorem{proposition}[subsection]{Proposition}

\theoremstyle{remark}
\newtheorem{example}[subsection]{Example}
\newtheorem{remark}[subsection]{Remark}

\usepackage{amsmath, amsfonts, amssymb, amsxtra}

\usepackage{ifpdf}
\ifpdf
\usepackage[pdftex]{graphics}
\IfFileExists{hyperref.sty}{\usepackage[pdftex]{hyperref}}{}
\else
\usepackage[dvips]{graphics}
\IfFileExists{hyperref.sty}{\usepackage[hypertex]{hyperref}}{}
\fi

\IfFileExists{euscript.sty}{\usepackage[mathcal]{euscript}}{}
\IfFileExists{mathrsfs.sty}{\usepackage{mathrsfs}}{\let\mathscr\mathfrak}

\newcommand{\ca}{{\mathcal A}}
\newcommand{\cb}{{\mathcal B}}
\newcommand{\cc}{{\mathcal C}}
\newcommand{\cd}{{\mathcal D}}
\newcommand{\ce}{{\mathcal E}}
\newcommand{\cm}{{\mathcal M}}
\newcommand{\cn}{{\mathcal N}}

\let\ge\geqslant
\let\le\leqslant

\let\tens\otimes
\let\eps\varepsilon

\newcommand{\tdt}{\otimes\dots\otimes}
\newcommand{\n}[1]{\nobreakdash-\hspace{0pt}}
\newcommand{\cQuiver}{{\mathscr Q}}
\newcommand{\ainf}[1]{$A_\infty$\nobreakdash-\hspace{0pt}}

\DeclareMathOperator\id{id}
\DeclareMathOperator\Ker{Ker}
\DeclareMathOperator\Ob{Ob}
\DeclareMathOperator\pr{pr}
\DeclareMathOperator\inj{in}

\begin{document}

\allowdisplaybreaks[1]

\title{Equalizers in the category of cocomplete cocategories}

\author{Bernhard Keller}
\address{UFR de Math\'ematiques\\
UMR 7586 du CNRS\\
Case 7012\\
Universit\'e Paris 7\\
2 place Jussieu\\
75251 Paris Cedex 05\\
France}
\email{keller@math.jussieu.fr}

\author{Oleksandr Manzyuk}
\address{Fachbereich Mathematik\\
Technische Universit\"at Kaiserslautern\\
Postfach 3049\\
Germany}
\email{manzyuk@mathematik.uni-kl.de}


\maketitle

\begin{abstract}
We prove existence of equalizers in certain categories of cocomplete
cocategories. This allows us to complete the proof
of the fact that \(A_\infty\)\n-functor categories arise as internal
Hom\n-objects in the category of differential graded cocomplete
augmented cocategories.
\end{abstract}

\section{Introduction}

We refer to \cite{math.RA/9910179} and to Lyubashenko-Ovsienko
\cite{LyuOvs-iResAiFn} for an introduction
to $A_\infty$-structures and their links to homological algebra.

The notion of an \(A_\infty\)\n-category appeared in Fukaya's work
on Floer homology \cite{Fukaya:A-infty}. Its relation to mirror
symmetry became apparent after Kontsevich's talk at ICM~'94
\cite{Kontsevich:alg-geom/9411018}. Following Kontsevich,
one should consider \(A_\infty\)\n-categories as models for
noncommutative varieties. This approach is being developed
by Kontsevich and Soibelman in \cite{math.RA/0606241}.

For a pair of \(A_\infty\)\n-categories \(\ca\) and \(\cb\), there is
an \(A_\infty\)\n-category \(A_\infty(\ca,\cb)\) whose objects are
\(A_\infty\)\n-functors and whose morphisms are
\(A_\infty\)\n-transformations. These \(A_\infty\)\n-functor categories
have been considered for example by Kontsevich \cite{Kontsevich:TriangCourse},
Fukaya \cite{Fukaya:FloerMirror-II}, Lef\`evre-Hasegawa
\cite{Lefevre-Ainfty-these}, Lyubashenko \cite{Lyu-AinfCat}. They
provide models for the internal Hom\n-functor of the homotopy category of
differential graded categories (Drinfeld \cite{Drinf:DGquot}, To\"en
\cite{math.AG/0408337}, cf.
\cite{math.KT/0601185} for a survey), where the internal Hom\n-functor is not a
derived functor. Furthermore, as detailed in \cite{math.RT/0510508},
\(A_\infty\)\n-functor categories yield a
natural construction of the \(B_\infty\)\n-structure on the Hochschild
complex of an associative algebra (Getzler-Jones
\cite{Getzler:hep-th/9403055}, Kadeishvili \cite{Kadeishvili88},
Voronov-Gerstenhaber \cite{GerstVoronov:Hochschild}), which is
important
for proving Deligne's conjecture and Tamarkin's version of
Kontsevich's formality theorem
(cf. for example Kontsevich-Soibelman \cite{math.QA/0001151}, Tamarkin
\cite{Tamarkin:math.QA/9803025}, Hinich \cite{Hinich03}).

In order to interpret \(A_\infty\)\n-functor categories as internal
Hom\n-objects, one passes to a suitable category of cocategories
following an idea of Lyubashenko \cite{Lyu-AinfCat}. For this suitable
category, one can either take the monoidal subcategory
generated by the images of graded quivers under the bar construction,
as in \cite{Lyu-AinfCat}, or the category of all (cocomplete augmented
etc.) cocategories. The former approach is developed further in the
forthcoming book by Bespalov, Lyubashenko, and Manzyuk
\cite{BesLyuMan-book} using the technique of closed multicategories. The
latter approach has been taken by the first author in
\cite{math.RT/0510508}. He proved in Theorem~5.3 of
\cite{math.RT/0510508} that the monoidal
category of cocomplete augmented cocategories was closed. However, the
proof of the theorem was incomplete: it relied on the assumption that
the category of cocomplete augmented cocategories has equalizers. In
this paper, we close this gap. Note that it suffices to prove existence
of equalizers in the category of cocomplete cocategories since it is
equivalent to the category of cocomplete augmented cocategories, see
Remark~\ref{rem-cocat-equiv-aug-cocat}.

\begin{theorem}
Suppose \(k\) is a field. Then the category of cocomplete
\(k\)\n-cocategories admits equalizers. The analogous assertions hold
in the graded and in the differential graded settings.
\end{theorem}

Since the proofs in the three cases are quite similar, we have chosen
to present the proof for the ungraded version providing remarks
concerning modifications necessary in the other cases. The proof
occupies Section~\ref{sec-proof}.

In the case of coalgebras, there is a different proof based on the
duality between coalgebras and algebras, successfully applied in works
of Kontsevich-Soibelman~\cite{math.RA/0606241} and
Hamilton-Lazarev~\cite{HamiltonLazarev}. We would like to briefly
outline it.

The category of finite dimensional coalgebras is anti-equivalent to the
category of finite dimensional algebras. Furthermore, finite
dimensional coalgebras are objects of finite presentation, in the
terminology of \cite[Definition~6.3.3]{KashiwaraSchapira05}, in the
category of coalgebras. Since an arbitrary coalgebra is a union of
finite dimensional subcoalgebras, see \cite{Green} or
\cite[Proposition~2.1.2]{math.RA/0606241}, it follows by
\cite[Proposition~6.3.4]{KashiwaraSchapira05} that the category of
coalgebras is equivalent to the category of ind-objects in the category
of finite dimensional coalgebras. Moreover, a finite dimensional
subcoalgebra of a cocomplete coalgebra is conilpotent, therefore the
category of cocomplete coalgebras is equivalent to the category of
ind-objects in the category of finite dimensional conilpotent
coalgebras, which is in turn anti-equivalent to the category of
pro-objects in the category of finite dimensional nilpotent algebras
(the category of formal algebras in the terminology of
\cite{HamiltonLazarev}). It suffices to establish the existence of
coequalizers in the latter category. However, by the dual of
\cite[Proposition~6.1.16]{KashiwaraSchapira05}, this follows from the
existence of coequalizers in the category of finite dimensional
nilpotent algebras.

Apparently, with some work the above argument can be generalized to
cocategories, although we did not check the details. The only drawback
of this approach, in our opinion, is that it is indirect. Our proof
relies on a direct verification and yields an explicit description of
equalizers, which is necessary in order to compute internal Hom-objects
in the category of cocomplete cocategories and to relate these to
\ainf-functor categories.

\bigskip\noindent
{\bf Acknowledgment.} The first author thanks V.~Hinich and the second
author for pointing out the gap which is at the origin of this note.

\section{Preliminaries}

Let \(k\) be a commutative ring. A \emph{\(k\)\n-quiver} \(\ca\) consists of a set of
objects \(\Ob\ca\) and of \(k\)\n-modules \(\ca(X,Y)\), for each pair
of objects \(X,Y\in\Ob\ca\). A \emph{morphism of \(k\)\n-quivers}
\(f:\ca\to\cb\) consists of a map \(\Ob f:\Ob\ca\to\Ob\cb\), \(X\mapsto
f(X)\), and of \(k\)\n-linear maps
\[
f=f_{X,Y}:\ca(X,Y)\to\cb(f(X),f(Y)),
\]
for each pair of objects \(X,Y\in\Ob\ca\). Let \(\cQuiver\) denote the
category of \(k\)\n-quivers. For a set \(S\), denote by \(\cQuiver/S\)
the subcategory of \(\cQuiver\) whose objects are \(k\)\n-quivers \(\ca\) such
that \(\Ob\ca=S\), and whose morphisms are morphisms of \(k\)\n-quivers
\(f \colon \ca\to\cb\) such that \(\Ob f=\id_S\). The category
\(\cQuiver/S\) is monoidal. The tensor product of quivers \(\ca\) and
\(\cb\) is given by
\[
(\ca\tens\cb)(X,Z)=\bigoplus_{Y\in S}\ca(X,Y)\tens\cb(Y,Z), \quad
X,Z\in S.
\]
The unit object is the \emph{discrete quiver} \(kS\) given by \(\Ob
kS=S\), \(kS(X,X)=k\) and \(kS(X,Y)=0\) if \(X\ne Y\), \(X,Y\in S\).
Recall that a \emph{cocategory} \((\cc,\Delta)\) is a coassociative
coalgebra in the monoidal category \(\cQuiver/\Ob\cc\). Thus, a
cocategory consists of a \(k\)\n-quiver \(\cc\) and of a morphism
\(\Delta\colon\cc\to\cc\tens\cc\) in \(\cQuiver/\Ob\cc\), the
comultiplication, satisfying the usual coassociativity condition. For
\(X,Y,Z\in\Ob\cc\), denote by
\[
\Delta_{X,Y,Z}=\pr_{X,Y,Z}\circ\Delta:\cc(X,Z)\to\cc(X,Y)\tens\cc(Y,Z)
\]
the components of \(\Delta\). Since the \(k\)\n-linear map
\[
(\Delta_{X,Y,Z})_{Y\in\Ob\cc}:\cc(X,Z)\to\prod_{Y\in\Ob\cc}\cc(X,Y)\tens\cc(Y,Z)
\]
factors as
\[
\cc(X,Z)\xrightarrow{\Delta}\bigoplus_{Y\in\Ob\cc}\cc(X,Y)\tens\cc(Y,Z)
\hookrightarrow\prod_{Y\in\Ob\cc}\cc(X,Y)\tens\cc(Y,Z),
\]
it follows that, for each \(t\in\cc(X,Z)\), the element \(\Delta_{X,Y,Z}(t)\) vanishes
for all but finitely many \(Y\in\Ob\cc\). The coassociativity is
expressed by the following equation:
\begin{multline*}
\bigl[
\cc(W,Z)\xrightarrow{\Delta_{W,X,Z}}\cc(W,X)\tens\cc(X,Z)\xrightarrow{1\tens\Delta_{X,Y,Z}}
\cc(W,X)\tens\cc(X,Y)\tens\cc(Y,Z)
\bigr]
\\
=\bigl[
\cc(W,Z)\xrightarrow{\Delta_{W,Y,Z}}\cc(W,Y)\tens\cc(Y,Z)\xrightarrow{\Delta_{W,X,Y}\tens1}
\cc(W,X)\tens\cc(X,Y)\tens\cc(Y,Z)
\bigr].
\end{multline*}
A \emph{cocategory homomorphism} \(f\colon(\cc,\Delta)\to(\cd,\Delta)\)
is a morphism of \(k\)\n-quivers \(f\colon\cc\to\cd\) compatible with
the comultiplication in the sense of the equation
\[
\bigl[
\cc\xrightarrow{f}\cd\xrightarrow{\Delta}\cd\tens\cd
\bigr]
=
\bigl[
\cc\xrightarrow{\Delta}\cc\tens\cc\xrightarrow{f\tens f}\cd\tens\cd
\bigr],
\]
where the morphism \(f\tens f\colon\cc\tens\cc\to\cd\tens\cd\) is given
by \(\Ob f\tens f=\Ob f\) and
\begin{multline*}
(f\tens f)_{X,Z}=\bigl[\bigoplus_{Y\in\Ob\cc}\cc(X,Y)\tens\cc(Y,Z)
\xrightarrow{\bigoplus_{Y\in\Ob\cc}f_{X,Y}\tens f_{Y,Z}}
\\
\bigoplus_{Y\in\Ob\cc}\cd(f(X),f(Y))\tens\cd(f(Y),f(Z))
\hookrightarrow\bigoplus_{U\in\Ob\cd}\cd(f(X),U)\tens\cd(U,f(Z))\bigr],
\end{multline*}
for each pair of objects \(X,Z\in\Ob\cc\). Explicitly, for \(X,Z\in\Ob\cc\),
\(U\in\Ob\cd\), the following equation holds true:
\begin{multline}
\bigl[
\cc(X,Z)\xrightarrow{f_{X,Z}}\cd(f(X),f(Z))\xrightarrow{\Delta_{f(X),U,f(Z)}}\cd(f(X),U)\tens\cd(U,f(Z))
\bigr]
\\
=\sum^{Y\in\Ob\cc}_{f(Y)=U}
\bigl[
\cc(X,Z)\xrightarrow{\Delta_{X,Y,Z}}\cc(X,Y)\tens\cc(Y,Z)\xrightarrow{f_{X,Y}\tens f_{Y,Z}}\cd(f(X),U)\tens\cd(U,f(Z))
\bigr].
\label{equ-f-Delta-Delta-f-f}
\end{multline}
In particular, the right hand side vanishes if \(U\) is not in the
image of \(f\).

Let \(\cc\) be a cocategory. Let \(\Delta^{(n)} \colon \cc\to\cc^{\tens n}\)
denote the comultiplication iterated \(n-1\) times, so that
\(\Delta^{(1)}=\id_\cc\), \(\Delta^{(2)}=\Delta\),
\(\Delta^{(3)}=(\Delta\tens1)\circ\Delta=(1\tens\Delta)\circ\Delta\), and so on.
Denote by
\[
\Delta^{(n)}_{X_0,\dots,X_n}=\pr_{X_0,\dots,X_n}\circ\Delta^{(n)} \colon \cc(X_0,X_n)\to\cc(X_0,X_1)\tens\cc(X_1,X_2)\tdt
\cc(X_{n-1},X_n)
\]
the components of \(\Delta^{(n)}\), for
\(X_0,\dots,X_n\in\Ob\cc\). Suppose \(f \colon \cc\to\cd\) is a cocategory
homomorphism. By induction on \(n\), it follows that
\begin{multline}
\Delta^{(n)}_{f(X),U_1,\dots,U_{n-1},f(Y)}\circ f_{X,Y}=\sum^{Z_1,\dots,Z_{n-1}\in\Ob\cc}_{f(Z_i)=U_i,\;i=1,\dots,n-1}
(f_{X,Z_1}\tens f_{Z_1,Z_2}\tdt
f_{Z_{n-1},Y})\circ\Delta^{(n)}_{X,Z_1,\dots,Z_{n-1},Y}:
\\
\cc(X,Y)\to\cd(f(X),U_1)\tens\cd(U_1,U_2)\tdt\cd(U_{n-1},f(Y)),
\label{equ-Delta-f-f-f-Delta}
\end{multline}
for an arbitrary collection of objects \(X,Y\in\Ob\cc\),
\(U_1,\dots,U_{n-1}\in\Ob\cd\).

A cocategory \(\cc\) is \emph{cocomplete}
if, for each pair of objects \(X,Y\in\Ob\cc\),
\[
\cc(X,Y)=\bigcup_{n\ge1}\Ker(\Delta^{(n)} \colon \cc(X,Y)\to\cc^{\tens
n}(X,Y)).
\]
Equivalently, \(\cc\) is cocomplete if for each \(t\in\cc(X,Y)\) there is \(n\ge1\) such that
\[
\Delta^{(n)}_{X,Z_1,\dots,Z_{n-1},Y}(t)=0,
\]
for all \(Z_1,\dots,Z_{n-1}\in\Ob\cc\).

\begin{example}
An arbitrary \(k\)\n-quiver \(\ca\) gives rise to a cocategory
\((T^{\ge1}\ca,\Delta)\), where \(T^{\ge1}\ca=\bigoplus_{n=1}^\infty
T^n\ca\), \(T^n\ca=\ca^{\tens n}\) is the \(n\)\n-fold tensor product
in \(\cQuiver/\Ob\ca\), and \(\Delta\) is the cut comultiplication.
Thus,
\[
T^{\ge1}\ca(X,Y)=\bigoplus_{Z_1,\dots,Z_{n-1}\in\Ob\ca}^{n\ge1}
\ca(X,Z_1)\tens\ca(Z_1,Z_2)\tdt\ca(Z_{n-1},Y),
\]
for each pair of objects \(X,Y\in\Ob\ca\), and \(\Delta\) is given by
\[
\Delta(f_1\tens f_2\tdt f_n)=\sum_{i=1}^{n-1}f_1\tdt
f_i\bigotimes f_{i+1}\tdt f_n.
\]
Since \(\Delta(T^n\ca)\subset\bigoplus_{p+q=n}^{p,q>0}T^p\ca\tens
T^q\ca\), it follows that  \(T^{\ge1}\ca\) is a cocomplete cocategory.
\end{example}

\begin{remark}
The correspondence \(\ca\mapsto T^{\ge1}\ca\) extends to a functor
\(T^{\ge1}:\cQuiver\to\cQuiver\). It is proven in
\cite[Chapter~8]{BesLyuMan-book} that the functor \(T^{\ge1}\) admits
the structure of a comonad, and that the category of
\(T^{\ge1}\)\n-coalgebras is isomorphic to the category of cocomplete
cocategories.
\end{remark}

\begin{remark}\label{rem-cocat-equiv-aug-cocat}
A cocategory \(\cc\) is \emph{counital} if it is equipped with a
morphism \(\eps:\cc\to k\Ob\cc\) in \(\cQuiver/\Ob\cc\) such that the
two counit equations hold. Note that, for an arbitrary set \(S\), the
\(k\)\n-quiver \(kS\) admits the natural structure of a counital
cocategory, namely the comultiplication is the canonical isomorphism
\(kS\xrightarrow{\sim}kS\tens kS\) in \(\cQuiver/S\) and the counit is
the identity map \(kS\to kS\). An \emph{augmented} cocategory is a
counital cocategory endowed with a morphism of counital cocategories
\(\eta:k\Ob\cc\to\cc\) such that \(\Ob\eta=\id_{\Ob\cc}\) and
\(\eps\circ\eta=\id_{k\Ob\cc}\). A morphism of augmented cocategories
is a cocategory homomorphism compatible with the counit and the
augmentation. The category of augmented cocategories is equivalent to
the category of cocategories \cite[Lemma~8.12]{BesLyuMan-book}:
given a cocategory \((\ca,\Delta)\), there is the natural structure of
an augmented cocategory on the \(k\)\n-quiver
\(T^{\le1}\ca=k\Ob\ca\oplus\ca\), where the counit and
the augmentation are the projection \(\eps=\pr_0:T^{\le1}\ca\to k\Ob\ca\)
and the inclusion \(\eta=\inj_0:k\Ob\ca\to T^{\le1}\ca\) respectively,
and the comultiplication is given by the formulas
\begin{align*}
\Delta|_{k\Ob\ca} =& \bigl[
k\Ob\ca\xrightarrow{\sim}
k\Ob\ca\tens k\Ob\ca\xrightarrow{\inj_0\tens\inj_0}T^{\le1}\ca\tens T^{\le1}\ca
\bigr],
\\
\Delta|_{\ca} =& \bigl[
\ca\xrightarrow{\sim}\ca\tens
k\Ob\ca\xrightarrow{\inj_1\tens\inj_0}T^{\le1}\ca\tens T^{\le1}\ca
\bigr]
\\
+& \bigl[
\ca\xrightarrow{\Delta}\ca\tens\ca\xrightarrow{\inj_1\tens\inj_1}T^{\le1}\ca\tens
T^{\le1}\ca
\bigr]
\\
+& \bigl[
\ca\xrightarrow{\sim}k\Ob\ca\tens
\ca\xrightarrow{\inj_0\tens\inj_1}T^{\le1}\ca\tens T^{\le1}\ca
\bigr].
\end{align*}
Conversely, given an augmented cocategory \((\cc,\Delta,\eps,\eta)\),
the reduced \(k\)\n-quiver \(\overline\cc=\Ker\eps\) becomes a
cocategory. The functors \(\ca\mapsto T^{\le1}\ca\) and
\(\cc\mapsto\overline\cc\) are quasi-inverse equivalences. By definition, an augmented
cocategory \(\cc\) is \emph{cocomplete} if its reduction
\(\overline\cc\) is cocomplete. Thus the category of cocomplete
cocategories is equivalent to the category of cocomplete augmented
cocategories.
\end{remark}

The above definitions admit obvious graded and differential graded
variants. For instance, a \emph{graded} (resp. \emph{differential
graded}) \emph{\(k\)\n-quiver} \(\ca\) consists of a set of objects
\(\Ob\ca\) and of graded \(k\)\n-modules (resp. cochain complexes of
\(k\)\n-modules) \(\ca(X,Y)\), for each pair of objects
\(X,Y\in\Ob\ca\). A \emph{morphism of graded} (resp. \emph{differential
graded}) \emph{\(k\)\n-quivers} \(f:\ca\to\cb\) consists of a map \(\Ob
f:\Ob\ca\to\Ob\cb\), \(X\mapsto f(X)\), and of morphisms of graded
\(k\)\n-modules of degree \(0\) (resp. cochain maps)
\[
f=f_{X,Y}:\ca(X,Y)\to\cb(f(X),f(Y)),
\]
for each pair of objects \(X,Y\in\Ob\ca\). The tensor product of
quivers with the same set of objects is defined analogously to the
considered case, using tensor product of graded \(k\)\n-modules (resp.
of cochain complexes). The definitions of cocategory, cocategory
homomorphism etc. are modified accordingly.

\section{Subcocategory generated by a set of objects}
\label{sec-subcocat}

Let \(\cc\) be a cocomplete cocategory, \(\cb\subset\cc\) a full subquiver, i.e.,
\(\Ob\cb\subset\Ob\cc\) and \(\cb(X,Y)=\cc(X,Y)\), for each pair of
objects \(X,Y\in\Ob\cb\). Then, in general, \(\cb\) is not a cocategory
since the comultiplication
\[
\Delta:\cc(X,Z)\to\bigoplus_{Y\in\Ob\cc}\cc(X,Y)\tens\cc(Y,Z)
\]
does not take values in \(\bigoplus_{Y\in\Ob\cb}\cc(X,Y)\tens\cc(Y,Z)\)
if \(X,Z\in\Ob\cb\). However, at least if \(k\) is a field, for each
subset \(S\subset\Ob\cc\) there exists a maximal cocomplete subcocategory
\(\cc_S\subset\cc\) such that \(\Ob\cc_S=S\). It is constructed as
follows. For a pair of objects \(X,Y\in S\), denote
\[
\cn(X,Y)=\bigoplus^{Z_1,\dots,Z_{n-1}\in\Ob\cc,\; n\ge1}_{\exists i \colon Z_i\not\in S}
\cc(X,Z_1)\tens\cc(Z_1,Z_2)\tdt\cc(Z_{n-1},Y),
\]
and define a \(k\)\n-linear map \(N=N_{X,Y} \colon \cc(X,Y)\to\cn(X,Y)\)
by
\[
N=\sum^{Z_1,\dots,Z_{n-1}\in\Ob\cc,\; n\ge1}_{\exists i \colon Z_i\not\in S}
\Delta^{(n)}_{X,Z_1,\dots,Z_{n-1},Y}.
\]
For each \(t\in\cc(X,Y)\), the sum in the right hand side is finite
since \(\cc\) is cocomplete. Let \(\cc_S(X,Y)=\Ker N\), so that we have
an exact sequence of \(k\)\n-vector spaces:
\[
0\to\cc_S(X,Y)\xrightarrow{\iota} \cc(X,Y)\xrightarrow{N}\cn(X,Y).
\]
Choose a splitting \(\pi \colon \cc(X,Y)\to\cc_S(X,Y)\) such that
\(\pi\circ\iota=\id_{\cc_S(X,Y)}\). Suppose that \(X,Y,Z\in S\). Then the composite
\[
\bigl[
\cc_S(X,Z)\xrightarrow{\iota}\cc(X,Z)\xrightarrow{\Delta_{X,Y,Z}}\cc(X,Y)\tens\cc(Y,Z)\xrightarrow{N\tens1}\cn(X,Y)\tens\cc(Y,Z)
\bigr]
\]
vanishes. Indeed, by coassociativity, \((N\tens1)\circ\Delta_{X,Y,Z}\) equals
\begin{multline*}
\sum^{Z_1,\dots,Z_{n-1}\in\Ob\cc,\; n\ge1}_{\exists i \colon Z_i\not\in S}
(\Delta^{(n)}_{X,Z_1,\dots,Z_{n-1},Y}\tens1)\circ\Delta_{X,Y,Z}
=\sum^{Z_1,\dots,Z_{n-1}\in\Ob\cc,\; n\ge1}_{\exists i \colon Z_i\not\in S}
\Delta^{(n+1)}_{X,Z_1,\dots,Z_{n-1},Y,Z}:
\\
\cc(X,Z)\to\bigoplus^{Z_1,\dots,Z_{n-1}\in\Ob\cc,\; n\ge1}_{\exists i \colon Z_i\not\in S}
\cc(X,Z_1)\tens\cc(Z_1,Z_2)\tdt\cc(Z_{n-1},Y)\tens\cc(Y,Z).
\end{multline*}
It follows from the definition of \(\iota\) that
\(\Delta^{(k)}_{X,U_1,\dots,U_{k-1},Y}\circ\iota=0\), for an arbitrary
sequence of objects \(U_1,\dots,U_{k-1}\in\Ob\cc\) such that
\(U_i\not\in S\) for some \(i\). Therefore
\((N\tens1)\circ\Delta_{X,Y,Z}\circ\iota=0\). Since the sequence
\[
0\to\cc_S(X,Y)\tens\cc(Y,Z)\xrightarrow{\iota\tens1}\cc(X,Y)\tens\cc(Y,Z)\xrightarrow{N\tens1}\cn(X,Y)\tens\cc(Y,Z)
\]
is exact, the composite \(\Delta_{X,Y,Z}\circ\iota\)
factors through \(\cc_S(X,Y)\tens\cc(Y,Z)\). In other words, there exists a unique
\(k\)\n-linear map \(\phi \colon \cc_S(X,Z)\to\cc_S(X,Y)\tens\cc(Y,Z)\) such that
\((\iota\tens1)\circ\phi=\Delta_{X,Y,Z}\circ\iota\). Since
\(\iota\tens1\) is an embedding split by \(\pi\tens1\), the map
\(\phi\) is necessarily given by the composite
\[
\bigl[
\cc_S(X,Z)\xrightarrow{\iota}\cc(X,Z)\xrightarrow{\Delta_{X,Y,Z}}\cc(X,Y)\tens\cc(Y,Z)\xrightarrow{\pi\tens1}\cc_S(X,Y)\tens\cc(Y,Z)
\bigr].
\]
In particular, the equation
\((\iota\tens1)\circ\phi=\Delta_{X,Y,Z}\circ\iota\) takes the form
\begin{equation}
(\iota\circ\pi\tens1)\circ\Delta_{X,Y,Z}\circ\iota=\Delta_{X,Y,Z}\circ\iota \colon \cc_S(X,Z)\to\cc(X,Y)\tens\cc(Y,Z).
\label{equ-iota-pi-1-Delta-iota-Delta-iota}
\end{equation}
Similarly, the following equation holds true:
\begin{equation}
(1\tens\iota\circ\pi)\circ\Delta_{X,Y,Z}\circ\iota=\Delta_{X,Y,Z}\circ\iota \colon \cc_S(X,Z)\to\cc(X,Y)\tens\cc(Y,Z).
\label{equ-1-iota-pi-Delta-iota-Delta-iota}
\end{equation}
Combining these equations yields
\[
(\iota\tens\iota)\circ(\pi\tens\pi)\circ\Delta_{X,Y,Z}\circ\iota=\Delta_{X,Y,Z}\circ\iota \colon \cc_S(X,Z)\to\cc(X,Y)\tens\cc(Y,Z).
\]
Define
\begin{equation}
\Delta'_{X,Y,Z}=\bigl[
\cc_S(X,Z)\xrightarrow{\iota}\cc(X,Z)\xrightarrow{\Delta_{X,Y,Z}}\cc(X,Y)\tens\cc(Y,Z)\xrightarrow{\pi\tens\pi}\cc_S(X,Y)\tens\cc_S(Y,Z)
\bigr].
\label{equ-Delta'}
\end{equation}
Then the above equation is equivalent to
\((\iota\tens\iota)\circ\Delta'_{X,Y,Z}=\Delta_{X,Y,Z}\circ\iota\). Coassociativity of
\(\Delta'\) follows from coassociativity of \(\Delta\) since
\begin{multline*}
(\iota\tens\iota\tens\iota)\circ(\Delta'\tens1)\circ\Delta'=(\Delta\tens1)\circ(\iota\tens\iota)\circ\Delta'=
(\Delta\tens1)\circ\Delta\circ\iota
\\
=(1\tens\Delta)\circ\Delta\circ\iota=(1\tens\Delta)\circ(\iota\tens\iota)\circ\Delta'=
(\iota\tens\iota\tens\iota)\circ(1\tens\Delta')\circ\Delta',
\end{multline*}
and \(\iota\tens\iota\tens\iota\) is an embedding split by
\(\pi\tens\pi\tens\pi\). Thus, \(\cc_S\) becomes a cocategory. The
quiver map \(\iota \colon \cc_S\to\cc\) with \(\Ob\iota \colon
S\hookrightarrow\Ob\cc\) is a cocategory homomorphism. Indeed, it was
shown above that equation~\eqref{equ-f-Delta-Delta-f-f} holds true for
\(X,Y,Z\in S\). If \(X,Z\in S\), \(Y\in\Ob\cc\smallsetminus S\), then
\(\Delta_{X,Y,Z}\circ\iota=0\) by the definition of \(\cc_S\),
therefore equation~\eqref{equ-f-Delta-Delta-f-f} is satisfied in this
case as well. By~\eqref{equ-Delta-f-f-f-Delta}, the equation
\(\iota^{\tens
n}\circ\Delta^{\prime(n)}_{X,Z_1,\dots,Z_{n-1},Y}=\Delta^{(n)}_{X,Z_1,\dots,Z_{n-1},Y}\circ\iota\)
holds true for all \(X,Z_1,\dots,Z_{n-1},Y\in S\). This implies that
the cocategory \(\cc_S\) is cocomplete: given an element
\(t\in\cc_S(X,Y)\), there is \(n\ge1\) such that
\(\Delta^{(n)}_{X,Z_1,\dots,Z_{n-1},Y}\circ\iota(t)=0\), for an
arbitrary collection of objects \(Z_1,\dots,Z_{n-1}\in\Ob\cc\). Since
\(\iota^{\tens n}\) is an embedding split by \(\pi^{\tens n}\), it
follows that \(\Delta^{\prime(n)}_{X,Z_1,\dots,Z_{n-1},Y}(t)=0\) for
all \(Z_1,\dots,Z_{n-1}\in S\).

\begin{proposition}[Universal property of \(\cc_S\)]
An arbitrary cocategory homomorphism \(h \colon \cb\to\cc\) with \(h(\Ob\cb)\subset
S\) factors uniquely through \(\cc_S\).
\end{proposition}

\begin{proof}
It follows from equation~\eqref{equ-Delta-f-f-f-Delta} that
\(\Delta^{(n)}_{h(X),U_1,\dots,U_{n-1},h(Y)}\circ h=0\) if \(U_i\not\in
S\) for some \(i\), thus \(N\circ h=0 \colon
\cb(X,Y)\to\cn(h(X),h(Y))\), for each pair \(X,Y\in\Ob\cb\). Therefore,
\(h \colon \cb(X,Y)\to\cc(h(X),h(Y))\) factors through
\(\cc_S(h(X),h(Y))\), i.e., there exists a unique linear map
\(\overline{h}=\overline{h}_{X,Y} \colon \cb(X,Y)\to\cc_S(h(X),h(Y))\)
such that \(\iota\circ\overline{h}=h\). The quiver map \(\overline{h}
\colon \cb\to\cc_S\) with \(\Ob\overline{h}=\Ob h \colon \Ob\cb\to S\)
is a cocategory homomorphism since
\[
(\iota\tens\iota)\circ\Delta'\circ\overline{h}=\Delta\circ\iota\circ\overline{h}=\Delta\circ
h=(h\tens
h)\circ\Delta=(\iota\tens\iota)\circ(\overline{h}\tens\overline{h})\circ\Delta,
\]
and \(\iota\tens\iota\) is an embedding split by \(\pi\tens\pi\).
Uniqueness of \(\overline{h}\) is obvious.
\end{proof}

\begin{remark}\label{rem-graded-diff-graded-case}
The same construction makes sense in graded and differential graded
contexts. The proofs transport literally except the following
subtlety in the case of differential graded cocategories: in general,
the embedding \(\iota:\cc_S(X,Y)\to\cc(X,Y)\) does not admit a
splitting which is a cochain map. Nevertheless, the argument can be
modified as follows. Choose a splitting \(\pi:\cc(X,Y)\to\cc_S(X,Y)\)
of \emph{graded} \(k\)\n-modules. Since \(\iota\) is a cochain map, i.e., \(d\circ\iota=\iota\circ
d\), it follows that the differential in \(\cc_S(X,Y)\) is necessarily
given by the composite
\[
d=\bigl[
\cc_S(X,Y)\xrightarrow{\iota}\cc(X,Y)\xrightarrow{d}\cc(X,Y)\xrightarrow{\pi}\cc_S(X,Y)
\bigr].
\]
In particular, the commutation relation \(d\circ\iota=\iota\circ d\)
takes the form
\[
\iota\circ\pi\circ d\circ\iota=d\circ\iota:\cc_S(X,Y)\to\cc(X,Y).
\]
Then the comultiplication \(\Delta'\) given by \eqref{equ-Delta'} is a
cochain map. Indeed,
\[
\Delta'\circ d=(\pi\tens\pi)\circ\Delta\circ\iota\circ\pi\circ
d\circ\iota=(\pi\tens\pi)\circ\Delta\circ d\circ\iota.
\]
On the other hand,
\begin{align*}
(1\tens d+d\tens1)\circ\Delta'=&(1\tens\pi\circ d\circ\iota+\pi\circ
d\circ\iota\tens1)\circ(\pi\tens\pi)\circ\Delta\circ\iota
\\
=&(\pi\tens\pi)\circ(1\tens d)\circ(1\tens\iota\circ\pi)\circ\Delta\circ\iota
\\
+&(\pi\tens\pi)\circ(d\tens1)\circ(\iota\circ\pi\tens1)\circ\Delta\circ\iota
\\
=&(\pi\tens\pi)\circ(1\tens d+d\tens1)\circ\Delta\circ\iota
\end{align*}
due to \eqref{equ-iota-pi-1-Delta-iota-Delta-iota} and
\eqref{equ-1-iota-pi-Delta-iota-Delta-iota}. Since \(\Delta\) is a
cochain map, it follows that \(\Delta'\circ d=(1\tens
d+d\tens1)\circ\Delta'\), thus \((\cc_S,\Delta')\) is a differential
graded cocategory. The further arguments remain unchanged.
\end{remark}

\section{Proof of the theorem}
\label{sec-proof}

Let \(\cc\), \(\cd\) be cocomplete cocategories, \(f,g \colon \cc\to\cd\)
cocategory homomorphisms. Denote \(S=\{X\in\Ob\cc|f(X)=g(X)\}\).
Suppose \(h \colon \cb\to\cc\) is a cocategory homomorphism such that \(f\circ
h=g\circ h\). Then for each \(W\in\Ob\cb\), \(f(h(W))=g(h(W))\),
therefore \(h(\Ob\cb)\subset S\). By the universal property of the
cocategory \(\cc_S\) there exists a unique cocategory homomorphism
\(\overline{h} \colon \cb\to\cc_S\) such that \(\iota\circ\overline{h}=h\).
Therefore, it suffices to construct an equalizer of the pair of
cocategory homomorphisms \(f\circ\iota,g\circ\iota \colon \cc_S\to\cd\). Thus,
we may assume without loss of generality that \(\Ob f=\Ob g\). Let us
construct an equalizer
\[
\ce\xrightarrow{e}{}\cc\overset{f}{\underset{g}{\rightrightarrows}}\cd
\]
in the category of cocomplete cocategories. Put \(\Ob\ce=\Ob\cc\). For
\(X,Y\in\Ob\cc\), denote by \(\cm(X,Y)\) the \(k\)\n-vector space
\[
\bigoplus_{\substack{X=X_0,X_1,\dots,X_p\in\Ob\cc\\ Y_0,\dots,Y_{q-1},Y_q=Y\in\Ob\cc}}^{p,q\ge0}
\cc(X,X_1)\tdt\cc(X_{p-1},X_p)\tens\cd(f(X_p),f(Y_0))\tens\cc(Y_0,Y_1)\tdt\cc(Y_{q-1},Y).
\]
Define a \(k\)\n-linear map \(R \colon \cc(X,Y)\to\cm(X,Y)\) by
\[
R=\sum_{\substack{X=X_0,X_1,\dots,X_p\in\Ob\cc\\ Y_0,\dots,Y_{q-1},Y_q=Y\in\Ob\cc}}^{p,q\ge0}
(1^{\tens p}\tens(f-g)\tens1^{\tens
q})\circ\Delta^{(p+1+q)}_{X,X_1,\dots,X_p,Y_0,\dots,Y_{q-1},Y}.
\]
It is well defined since \(\cc\) is cocomplete. Let \(\ce(X,Y)=\Ker R\), so that we have an exact sequence
\[
0\to\ce(X,Y)\xrightarrow{e}\cc(X,Y)\xrightarrow{R} \cm(X,Y).
\]
Choose a splitting \(p \colon \cc(X,Y)\to\ce(X,Y)\) such that \(p\circ
e=\id_{\ce(X,Y)}\). Suppose that \(X,Y,Z\in\Ob\cc\). Then the composite
\[
\bigl[\ce(X,Z)\xrightarrow{e}\cc(X,Z)\xrightarrow{\Delta_{X,Y,Z}}\cc(X,Y)\tens\cc(Y,Z)\xrightarrow{R\tens1}\cm(X,Y)\tens\cc(Y,Z)\bigr]
\]
vanishes. Indeed, by coassociativity, the composite
\((R\tens1)\circ\Delta_{X,Y,Z}\) equals
\begin{multline*}
\sum_{\substack{X=X_0,X_1,\dots,X_p\in\Ob\cc\\ Y_0,\dots,Y_{q-1},Y_q=Y\in\Ob\cc}}^{p,q\ge0}
(1^{\tens p}\tens(f-g)\tens1^{\tens
q+1})\circ(\Delta^{(p+1+q)}_{X,X_1,\dots,X_p,Y_0,\dots,Y_{q-1},Y}\tens1)\circ\Delta_{X,Y,Z}
\\
\hfill=\sum_{\substack{X=X_0,X_1,\dots,X_p\in\Ob\cc\\ Y_0,\dots,Y_{q-1},Y_q=Y\in\Ob\cc}}^{p,q\ge0}
(1^{\tens p}\tens(f-g)\tens1^{\tens
q+1})\circ\Delta^{(p+1+q+1)}_{X,X_1,\dots,X_p,Y_0,\dots,Y_{q-1},Y,Z} \colon
\hskip\multlinegap
\\
\hskip\multlinegap\cc(X,Z)\to\bigoplus_{\substack{X=X_0,X_1,\dots,X_p\in\Ob\cc\\ Y_0,\dots,Y_{q-1},Y_q=Y\in\Ob\cc}}^{p,q\ge0}
\cc(X,X_1)\tdt\cc(X_{p-1},X_p)\tens\cd(f(X_p),f(Y_0))\hfill
\\
\tens\cc(Y_0,Y_1)\tdt\cc(Y_{q-1},Y)\tens\cc(Y,Z).
\end{multline*}
It follows from the definition \(\ce\) that \((1^{\tens
p}\tens(f-g)\tens1^{\tens
q+1})\circ\Delta^{(p+1+q+1)}_{X,X_1,\dots,X_p,Y_0,\dots,Y_{q-1},Y,Z}\circ
e=0\), for all \(X,X_1,\dots,X_p,Y_0,\dots,Y_{q-1},Y,Z\in\Ob\cc\),
therefore \((R\tens1)\circ\Delta_{X,Y,Z}\circ e=0\). Since the sequence
\[
0\to\ce(X,Y)\tens\cc(Y,Z)\xrightarrow{e\tens1}\cc(X,Y)\tens\cc(Y,Z)\xrightarrow{R\tens1}\cm(X,Y)\tens\cc(Y,Z)
\]
is exact, it follows that the map \(\Delta_{X,Y,Z}\circ e\)
factors through \(\ce(X,Y)\tens\cc(Y,Z)\). In other words, there exists a unique
\(k\)\n-linear map \(\psi \colon \ce(X,Z)\to\ce(X,Y)\tens\cc(Y,Z)\) such that
\((e\tens1)\circ\psi=\Delta_{X,Y,Z}\circ e\). Since
\(e\tens1\) is an embedding split by \(p\tens1\), the map
\(\psi\) is necessarily given by the composite
\[
\bigl[
\ce(X,Z)\xrightarrow{e}\cc(X,Z)\xrightarrow{\Delta_{X,Y,Z}}\cc(X,Y)\tens\cc(Y,Z)\xrightarrow{p\tens1}\ce(X,Y)\tens\cc(Y,Z)
\bigr].
\]
In particular, the equation
\((e\tens1)\circ\psi=\Delta_{X,Y,Z}\circ e\) takes the form
\[
(e\circ p\tens1)\circ\Delta_{X,Y,Z}\circ e=\Delta_{X,Y,Z}\circ e \colon \ce(X,Z)\to\cc(X,Y)\tens\cc(Y,Z).
\]
Similarly, the following equation holds true:
\[
(1\tens e\circ p)\circ\Delta_{X,Y,Z}\circ e=\Delta_{X,Y,Z}\circ e \colon \ce(X,Z)\to\cc(X,Y)\tens\cc(Y,Z).
\]
Combining these equations yields
\[
(e\tens e)\circ(p\tens p)\circ\Delta_{X,Y,Z}\circ e=\Delta_{X,Y,Z}\circ e \colon \ce(X,Z)\to\cc(X,Y)\tens\cc(Y,Z).
\]
Define
\[
\tilde\Delta_{X,Y,Z}=\bigl[
\ce(X,Z)\xrightarrow{e}\cc(X,Z)\xrightarrow{\Delta_{X,Y,Z}}\cc(X,Y)\tens\cc(Y,Z)\xrightarrow{p\tens p}\ce(X,Y)\tens\ce(Y,Z)
\bigr].
\]
As in the case of \(\cc_S\), one shows that \(\tilde\Delta\) turns
\(\ce\) into a cocomplete cocategory, and that \(e\) becomes a cocategory
homomorphism.

Suppose \(h \colon \cb\to\cc\) is a cocategory homomorphism such that
\(f\circ h=g\circ h\). Then \(R\circ h=0 \colon \cb(X,Y)\to\cm(h(X),h(Y))\).
Indeed, by identity~\eqref{equ-Delta-f-f-f-Delta},
\begin{multline*}
\Delta^{(p+1+q)}_{h(X),X_1,\dots,X_p,Y_0,\dots,Y_{q-1},h(Y)}\circ h
\\
=\sum^{h(U_1)=X_1,\dots,h(U_p)=X_p}_{h(V_0)=Y_0,\dots,h(V_{q-1})=Y_{q-1}}(h_{X,U_1}\tdt h_{U_p,V_0}\tdt h_{V_{q-1},Y})
\circ\Delta^{(p+1+q)}_{X,U_1,\dots,U_p,V_0,\dots,V_{q-1},Y}.
\end{multline*}
It follows that
\[
R\circ h=
\sum^{h(U_1)=X_1,\dots,h(U_p)=X_p}_{h(V_0)=Y_0,\dots,h(V_{q-1})=Y_{q-1}}
(h^{\tens p}\tens (f-g)\circ h\tens h^{\tens
q})\circ\Delta^{(p+1+q)}_{X,U_1,\dots,U_p,V_0,\dots,V_{q-1},Y}=0,
\]
for \(X_1,\dots,X_p,Y_0,\dots,Y_{q-1}\in\Ob\cc\). Therefore the map
\(h \colon \cb(X,Y)\to\cc(h(X),h(Y))\) factors through
\(e \colon \ce(h(X),h(Y))\to\cc(h(X),h(Y))\), i.e., \(h=e\circ j\) for some
\(k\)\n-linear map \(j \colon \cb(X,Y)\to\ce(h(X),h(Y))\). The morphism of \(k\)\n-quivers \(j \colon \cb\to\ce\) with
\(\Ob j=\Ob h\) is a cocategory homomorphisms since
\[
(e\tens
e)\circ(j\tens j)\circ\Delta=(h\tens h)\circ\Delta =\Delta\circ
h=\Delta\circ e\circ j=(e\tens e)\circ\Delta\circ j,
\]
and \(e\tens e\) is an embedding split by \(p\tens p\). Uniqueness of \(j\) is obvious since
\(e\) is an embedding.\qed

\begin{remark}
The theorem is true for cocomplete graded (resp. differential graded)
cocategories as well, with appropriate modifications in the proof similar to
those made in Remark~\ref{rem-graded-diff-graded-case}.
\end{remark}

\begin{remark}
The same proof shows the existence of equalizers in the category of
cocomplete coalgebras, which are just cocategories with only one
object. The intermediate step described in Section~\ref{sec-subcocat}
becomes superfluous.
\end{remark}

\providecommand{\bysame}{\leavevmode\hbox to3em{\hrulefill}\thinspace}
\providecommand{\MR}{\relax\ifhmode\unskip\space\fi MR }
\providecommand{\MRhref}[2]{%
  \href{http://www.ams.org/mathscinet-getitem?mr=#1}{#2}
}
\providecommand{\href}[2]{#2}


\begin{thebibliography}{10}

\bibitem{BesLyuMan-book}
Yuri Bespalov, V.~V. Lyubashenko, and Oleksandr Manzyuk,
\emph{Closed multicategory of pretriangulated ${A}_\infty$-categories},
book in progress, 2006, \texttt{http://www.math.ksu.edu/\(\sim\)lub/papers.html}.

\bibitem{Drinf:DGquot}
Vladimir~G. Drinfeld, \emph{{DG} quotients of {DG} categories}, J. Algebra
  \textbf{272} (2004), no.~2, 643--691,
  \texttt{math.KT/\linebreak[1]0210114}.

\bibitem{Fukaya:A-infty}
Kenji Fukaya, \emph{Morse homotopy, ${A}_\infty$-category, and {F}loer
  homologies}, {Proc. of GARC Workshop on Geometry and Topology '93} (H.~J.
  Kim, ed.), Lecture Notes, no.~18, Seoul Nat. Univ., Seoul, 1993,
  pp.~1--102, \texttt{http://www.math.kyoto-u.ac.jp/\(\sim\)fukaya/fukaya.html}.

\bibitem{Fukaya:FloerMirror-II}
\bysame, \emph{Floer homology and mirror symmetry. {II}}, {Minimal surfaces,
  geometric analysis and symplectic geometry (Baltimore, MD, 1999)}, Adv. Stud.
  Pure Math., vol.~34, Math. Soc. Japan, Tokyo, 2002, pp.~31--127.

\bibitem{Getzler:hep-th/9403055}
Ezra Getzler and John D.~S. Jones, \emph{Operads, homotopy algebra, and
  iterated integrals for double loop spaces}, 1994,
  \texttt{hep-th/\linebreak[1]9403055}.

\bibitem{Green}
J.~A. Green, \emph{Locally finite representations}, J. Algebra
\textbf{41} (1976), 137--171.

\bibitem{HamiltonLazarev}
Alastair Hamilton and Andrey Lazarev, \emph{Homotopy algebras and
noncommutative geometry}, 2004, \texttt{math.QA/\linebreak[1]0410621}.

\bibitem{Hinich03}
Vladimir Hinich, \emph{Tamarkin's proof of {K}ontsevich formality theorem},
  Forum Math. \textbf{15} (2003), no.~4, 591--614,
  \texttt{math.QA/\linebreak[1]0003052}.

\bibitem{Kadeishvili88}
Tornike~V. Kadeishvili, \emph{The structure of the ${A}(\infty)$-algebra, and
  the {H}ochschild and {H}arrison cohomologies}, Proc. of A. Razmadze Math.
  Inst. \textbf{91} (1988), 20--27,
  \texttt{math.AT/\linebreak[1]0210331}.

\bibitem{KashiwaraSchapira05}
Masaki Kashiwara and Pierre Schapira, \emph{Categories and sheaves},
  Grundlehren der Mathematischen Wissenschaften [Fundamental Principles of
  Mathematical Sciences], vol. 332, Springer-Verlag, Berlin, 2005.

\bibitem{math.RA/9910179}
Bernhard Keller, \emph{Introduction to {A}-infinity algebras and modules},
  Homology, Homotopy and Applications \textbf{3} (2001), no.~1, 1--35,
  \texttt{math.RA/\linebreak[1]9910179},
  \texttt{http://intlpress.com/HHA/v3/n1/a1/}.

\bibitem{math.RT/0510508}
\bysame, \emph{A-infinity algebras, modules and functor categories}, 2005,
\texttt{math.RT/\linebreak[1]0510508}.

\bibitem{math.KT/0601185}
\bysame, \emph{On differential graded categories}, contribution to the Proceedings
of the ICM 2006, \texttt{math.KT/\linebreak[1]0601185}.

\bibitem{Kontsevich:alg-geom/9411018}
Maxim Kontsevich, \emph{Homological algebra of mirror symmetry}, {Proc.
  Internat. Cong. Math., Z\"urich, Switzerland 1994} (Basel), vol.~1,
  Birkh\"auser Verlag, 1995, pp.~120--139, \texttt{math.AG/\linebreak[1]9411018}.

\bibitem{Kontsevich:TriangCourse}
\bysame, \emph{Triangulated categories and geometry},
  Course at the \'Ecole Normale Sup\'eriure, Paris, March and April
1998, available at \texttt{http://www.math.uchicago.edu/\(\sim\)arinkin/}.

\bibitem{math.QA/0001151}
Maxim Kontsevich and Yan~S. Soibelman, \emph{Deformations of algebras over
  operads and {D}eligne's conjecture}, {Conf\'erence Mosh\'e Flato 1999, Vol. I
  (Dijon)}, Math. Phys. Stud., vol.~21, Kluwer Academic Publishers, Dordrecht,
  2000, pp.~255--307, \texttt{math.QA/\linebreak[1]0001151}.

\bibitem{math.RA/0606241}
\bysame, \emph{Notes on {A}-infinity algebras, {A}-infinity categories and
  non-commutative geometry. {I}}, 2006,
  \texttt{math.RA/\linebreak[1]0606241}.

\bibitem{Lefevre-Ainfty-these}
Kenji Lef\`evre-Hasegawa, \emph{Sur les ${A}_\infty$-cat\'egories}, Ph.D.
  thesis, Universit\'e Paris 7, U.F.R. de Math\'ematiques, 2003,
  \texttt{math.CT/\linebreak[1]0310337}.

\bibitem{Lyu-AinfCat}
V.~V. Lyubashenko, \emph{Category of ${A}_\infty$-categories}, Homology,
  Homotopy and Applications \textbf{5} (2003), no.~1, 1--48,
  \texttt{math.CT/\linebreak[1]0210047},
  \texttt{http://intlpress.com/HHA/v5/n1/a1/}.

\bibitem{LyuOvs-iResAiFn}
V.~V. Lyubashenko and Serge~A. Ovsienko, \emph{A construction of quotient
  ${A}_\infty$-categories}, Homology, Homotopy and Applications \textbf{8}
  (2006), no.~2, 157--203, \texttt{math.CT/\linebreak[1]0211037},
  \texttt{http://intlpress.com/HHA/v8/n2/a9/}.

\bibitem{Tamarkin:math.QA/9803025}
Dmitry~E. Tamarkin, \emph{Another proof of {M}.~{K}ontsevich formality
  theorem}, 1998, \texttt{math.QA/\linebreak[1]9803025}.

\bibitem{math.AG/0408337}
Bertrand To\"en, \emph{The homotopy theory of dg-categories and derived Morita
  theory}, 2004, \texttt{math.AG/\linebreak[1]0408337}.

\bibitem{GerstVoronov:Hochschild}
Alexander~A. Voronov and Murray Gerstenhaber, \emph{Higher operations on the
  {H}ochschild complex}, Functional Anal. Appl. \textbf{29} (1995), no.~1,
  1--6.

\end{thebibliography}

\end{document}